\documentclass[11pt]{article}

\usepackage{amsmath, amsthm, amssymb}
\usepackage{enumitem}
\usepackage{pdflscape}
\usepackage{caption}
\usepackage{bm}
\usepackage{ifpdf}
	\ifpdf
\usepackage[pdftex]{graphicx}
	\else
\usepackage[dvips]{graphicx}
	\fi
\usepackage[all]{xy}
\usepackage{tocvsec2}
\usepackage{bbm}
	\input xy
	\xyoption{all}
\usepackage[pdftex,plainpages=false,hypertexnames=false,pdfpagelabels]{hyperref}
	\newcommand{\arxiv}[1]{\href{http://arxiv.org/abs/#1}{\tt arXiv:\nolinkurl{#1}}}
	\newcommand{\arXiv}[1]{\href{http://arxiv.org/abs/#1}{\tt arXiv:\nolinkurl{#1}}}
	\newcommand{\doi}[1]{\href{http://dx.doi.org/#1}{{\tt DOI:#1}}}

	\newcommand{\googlebooks}[1]{(preview at \href{http://books.google.com/books?id=#1}{google books})}
	
\usepackage{xcolor}
	\definecolor{dark-red}{rgb}{0.7,0.25,0.25}
	\definecolor{dark-blue}{rgb}{0.15,0.15,0.55}
	\definecolor{medium-blue}{rgb}{0,0,.8}
	\definecolor{DarkGreen}{RGB}{0,150,0}
	\definecolor{rho}{named}{red}
	\hypersetup{
	   colorlinks, linkcolor={purple},
	   citecolor={medium-blue}, urlcolor={medium-blue}
	}
\usepackage{longtable}
\usepackage{fullpage}

\usepackage{mathbbol}

\theoremstyle{plain}
\newtheorem{thm}{Theorem}[section]
\newtheorem*{thm*}{Theorem}
\newtheorem{thmalpha}{Theorem}

\newtheorem{cor}[thm]{Corollary}

\newtheorem*{cor*}{Corollary}

\newtheorem*{conj*}{Conjecture}
\newtheorem{lem}[thm]{Lemma}
\newtheorem{prop}[thm]{Proposition}

\newtheorem*{quest*}{Question}
\newtheorem*{claim*}{Claim}

\theoremstyle{definition}

\newtheorem{construction}[thm]{Construction}

\newtheorem{exs}[thm]{Examples}

\newtheorem*{ex*}{Example}
\newtheorem{sub-ex}[thm]{Sub-Example}
\newtheorem{counter-ex}[thm]{Counter-Example}

\newtheorem*{rem*}{Remark}
\newtheorem{remark}[thm]{Remark}

\DeclareMathOperator{\Ad}{Ad}
\DeclareMathOperator{\Aut}{Aut}

\DeclareMathOperator{\End}{End}
\DeclareMathOperator{\ev}{ev}
\DeclareMathOperator{\Hom}{Hom}

\DeclareMathOperator{\op}{op}

\DeclareMathOperator{\id}{id}

\newcommand{\comment}[1]{}

\newcommand{\noshow}[1]{}
\newcommand{\MR}[1]{}

\newcommand{\fgpBim}{{\sf Bim_{fgp}}}
\newcommand{\fgpMod}{{\sf Mod_{fgp}}}

\renewcommand{\Vec}{{\sf Vec}}

\newcommand{\fdHilb}{{\sf Hilb_{fd}}}
\newcommand{\rCorr}{{\mathsf{C^{*}Alg}}}

\newcommand{\Ind}{{\sf Ind}}

\newcommand{\CDisc}{{\sf C^*Disc}}

\newcommand{\PQR}{{\sf PQR}}
\newcommand{\PQN}{{\sf PQN}}

\renewcommand{\>}{\rangle}

\def\semicolon{;}
\def\applytolist#1{
    \expandafter\def\csname multi#1\endcsname##1{
        \def\multiack{##1}\ifx\multiack\semicolon
            \def\next{\relax}
        \else
            \csname #1\endcsname{##1}
            \def\next{\csname multi#1\endcsname}
        \fi
        \next}
    \csname multi#1\endcsname}

\def\calc#1{\expandafter\def\csname c#1\endcsname{{\mathcal #1}}}
\applytolist{calc}QWERTYUIOPLKJHGFDSAZXCVBNM;
\def\bbc#1{\expandafter\def\csname bb#1\endcsname{{\mathbb #1}}}
\applytolist{bbc}QWERTYUIOPLKJHGFDSAZXCVBNM;
\def\bfc#1{\expandafter\def\csname bf#1\endcsname{{\mathbf #1}}}
\applytolist{bfc}QWERTYUIOPLKJHGFDSAZXCVBNM;
\def\sfc#1{\expandafter\def\csname s#1\endcsname{{\sf #1}}}
\applytolist{sfc}QWERTYUIOPLKJHGFDSAZXCVBNM;
\def\fc#1{\expandafter\def\csname f#1\endcsname{{\mathfrak #1}}}
\applytolist{fc}QWERTYUIOPLKJHGFDSAZXCVBNM;

\usepackage{tikz}
\usepackage{tikz-cd}
\usetikzlibrary{arrows,backgrounds,patterns.meta}
\usetikzlibrary{positioning,shadings,cd}
\usetikzlibrary{shapes}
\usetikzlibrary{backgrounds}
\usetikzlibrary{decorations,decorations.pathreplacing,decorations.markings}
\usetikzlibrary{fit,calc,through}
\usetikzlibrary{external}
\usetikzlibrary{arrows}
\tikzset{vertex/.style = {shape=circle,draw,fill=black,inner sep=0pt,minimum size=5pt}}
\tikzset{edge/.style = {->,> = latex', bend right}}
\tikzset{
	super thick/.style={line width=3pt}
}
\tikzset{
    quadruple/.style args={[#1] in [#2] in [#3] in [#4]}{
        #1,preaction={preaction={preaction={draw,#4},draw,#3}, draw,#2}
    }
}
\tikzstyle{shaded}=[fill=red!10!blue!20!gray!30!white]
\tikzstyle{unshaded}=[fill=white]
\tikzstyle{empty box}=[circle, draw, thick, fill=white, opaque, inner sep=2mm]
\tikzstyle{annular}=[scale=.7, inner sep=1mm, baseline]
\tikzstyle{rectangular}=[scale=.75, inner sep=1mm, baseline=-.1cm]
\tikzstyle{mid>}=[decoration={markings, mark=at position 0.5 with {\arrow{>}}}, postaction={decorate}]
\tikzstyle{mid<}=[decoration={markings, mark=at position 0.5 with {\arrow{<}}}, postaction={decorate}]
\tikzstyle{over}=[double, draw=white, super thick, double=]

\tikzdeclarepattern{
  name=primeddots,
  type=uncolored,
  bounding box={(-.6pt,-.6pt) and (.6pt,.6pt)},
  tile size={(5pt,5pt)},
  tile transformation={rotate=60},
  parameters={none}, 
  code={
    \fill(0pt,0pt) circle (.35pt);
  }
}

\tikzdeclarepattern{
  name=primeddots2,
  type=uncolored,
  bounding box={(-.3pt,-.3pt) and (.3pt,.3pt)},
  tile size={(2.5pt,2.5pt)},
  tile transformation={rotate=60},
  parameters={none}, 
  code={
    \fill(0pt,0pt) circle (.35pt);
  }
}

\tikzstyle{primedregion}[none]=[
	preaction={fill=#1},
	pattern=primeddots,
  draw=#1,
]

\tikzstyle{primedregion2}[none]=[
	preaction={fill=#1},
	pattern=primeddots2,
  draw=#1,
]

\tikzcdset{scale cd/.style={every label/.append style={scale=#1},
    cells={nodes={scale=#1}}}}

\let\OLDthebibliography\thebibliography
\renewcommand\thebibliography[1]{
  \OLDthebibliography{#1}
  \setlength{\parskip}{0pt}
  \setlength{\itemsep}{0pt plus 0.3ex}
}

\begin{document}

\title{Discrete inclusions from Cuntz-Pimsner algebras}
\author{Roberto Hern\'{a}ndez Palomares}
\date{\today}
\maketitle
\begin{abstract}
We show that the core inclusion arising from a Cuntz-Pimsner algebra generated by a full, faithful and dualizable correspondence is C*-discrete, and express it as a crossed-product by an action of a unitary tensor category.
In particular, we show $\mathsf{UHF}_{n^\infty}\subset \mathcal{O}_n$, the inclusion of the fixed-point algebra under the gauge symmetry of the Cuntz algebra, is irreducible and C*-discrete.
We describe the dualizable bimodules appearing under this inclusion, including their semisimple decompositions and fusion rules, their Watatani indices and Pimsner-Popa bases, as well as their sets of cyclic algebraic generators. 
\end{abstract}

\section{Introduction}

Taking inspiration from subfactors, to a C*-correspondence $X$ over a unital C*-algebra $D$, Pimsner associated the so-called Cuntz-Pimsner algebra $\cO_X$ that simultaneously generalizes Cuntz-Krieger algebras and crossed-products by single automorphisms \cite{MR1426840}. 
Later, Exel showed that Cuntz-Krieger algebras determined by $n\times n$ nondegenerate matrices arise from a Markov subshift acting on a compact set via an endomorphic crossed-product construction \cite[Theorem 6.2]{MR2032486}.  
Our purpose here is to incorporate various crossed-product constructions by non-invertible dynamics arising from Cuntz-Pimsner algebras into the framework of unitary tensor categories (UTC) and their actions on C*-algebras. 
In doing so, we provide more explicit examples of actions of UTCs on C*-algebras, which are not widely available in the literature.  
Indeed, the crossed-products by single imprimitivity bimodules \cite{MR1467459, MR2029622} or by single unital endomorphisms \cite{MR2032486, MR1966826} often organically fall into the broader framework of \emph{C*-discrete inclusions} given by crossed-products by tensor categories acting by \emph{dualizable}/\emph{finitely generated projective} bimodules.

A unital inclusion of C*-algebras $A\overset{E}{\subset}B$, where $E:B\twoheadrightarrow A$ is a faithful conditional expectation is C*-discrete if $B$ is densely spanned by a family of finitely generated projective $A$-$A$ bimodules.
It was established in \cite[Theorem A]{2023arXiv230505072H} that $A\subset B$ is \emph{irreducible} (i.e. $A'\cap B\cong \bbC1$) and C*-discrete if and only if there is an \emph{outer action} $F:\cC\to \fgpBim(A)$ of $\cC$ on $A$, where $\cC$ is a UTC and $F$ a fully-faithful unitary tensor functor, a $\cC$-graded C*-algebra object $\bbB$ which is \emph{connected} (i.e., $\bbB(1_\cC)\cong \bbC$),
and there exists an expectation-preserving $*$-isomorphism 
$
(A\subset B)\ \cong\ (A\subset A\rtimes_{r,F}\bbB).
$  
Here, the expectation on the crossed-product projects  onto the $1_\cC$-graded component.

Our main result establishes that if $X$ is a full, faithful and  dualizable bimodule over a unital C*-algebra $D$ (e.g. if $X={}_{N}L^2(M)_N$ for a finite-index subfactor $N\subset M$ \cite{MR696688}), then the ``extended scalars'' $\varinjlim\End(X^{\boxtimes n}_D)=\cO_X^\bbT$ roughly behave like subfactors of $\cO_X.$
More precisely, the associated inclusion arising as fixed points over the circle's gauge action is C*-discrete.
\begin{thmalpha}[{Theorem~\ref{thm:CuntzPimsnerCoreDisc}, Corollary~\ref{cor:CPCrossedProduct}}]\label{thmalpha:CuntzPimsnerCoreDisc}
    If $D$ is a unital C*-algebra, and $X\in \fgpBim(D)$ is full with unital (hence isometric) left $D$-action, then $\cO_X^\bbT \overset{E}{\subset} \cO_X$ is a C*-discrete inclusion. 
    Moreover, 
    $$
    \cO_X\cong \cO_X^\bbT\rtimes_{r,F}\bbC[\bbZ]
    $$
    is a reduced crossed-product by a (not necessarily outer) action $F$ of $\bbZ$ on $\cO_X^\bbT$ by dualizable bimodules under the group C*-algebra object $\bbC[\bbZ].$
\end{thmalpha}
\noindent Being next in complexity to a crossed-product by a single endomorphism, the inclusions  $\cO^\bbT_X\subset \cO_X$ are natural targets to study structural properties preserved under crossed-products by a single dualizable correspondence (c.f. \cite[\S7]{MR2102572}, \cite[Theorem 4.2]{2024arXiv240918161H}). 
We highlight that Cuntz-Krieger algebras arising from $n\times n$ matrices with no nonzero rows and columns fall under the scope of this theorem (c.f. Examples \ref{exs:C,CK}).

To prove Theorem \ref{thmalpha:CuntzPimsnerCoreDisc} we rely on the $\bbZ$-graded $*$-isomorphism $\cO_X \cong \cO_\rho$ constructed in \cite{MR1658088}. where $\cO_\rho$ denotes the Doplicher-Roberts algebra \cite{MR1010160} generated by the object $\rho = X$ considered as a right Hilbert $D$-module. 
This isomorphism takes the $k$-th spectral subspace of the gauge action on $\cO_X$ onto the $k$-th graded component of $\cO_\rho$ given by an inductive limit of certain fgp $D$-$D$ bimodules. 
We then show that the limit bimodule remains dualizable when considered over the extended scalars, thus witnessing the \emph{projective quasi-normalizer} \cite[Definition 2.3]{2023arXiv230505072H} of the inclusion is dense in $\cO_\rho.$ 
One can also view $\cO_\rho$ as a compression of the Cuntz-Krieger graph algebra of the (bidirected) fusion graph of $X$ \cite{MR3624399}, but we will not rely on this presentation.

In this generality, it is not necessarily true that the gauge inclusion is irreducible, nor that $\cO_X^\bbT$ has trivial center. 
Moreover, since $D$ is not assumed to have trivial center, the associated categories of bimodules will not be UTC's but rather rigid C*-tensor categories whose units may have nontrivial endomorphisms, and providing a detailed description of its bimodules may be hard.

\medskip

Assuming  $D=\bbC$ and $X=\bbC^n$ we recover the Cuntz algebras $\cO_n$ for $n\in \bbN,$ and in this case we can show that the inclusions from Cuntz algebra cores by their canonical gauge action are irreducible and C*-discrete. And thus, Cuntz algebras are crossed-products by an action of a tensor category: 
\begin{thmalpha}[{Theorem~\ref{thm:CuntzCoreDisc}, Corollary~\ref{cor:OnCrossedProduct}}]\label{thmalpha:CuntzCoreDisc}
    For any $2\leq n\in \bbN,$ the inclusion $\cO_n^\bbT\overset{E}{\subset}\cO_n$ arising from the fixed points of the canonical gauge action is irreducible and C*-discrete. 
    Moreover, 
    $$\cO_n\cong \cO_n^\bbT\rtimes_{r,F} \bbC[\bbZ]
    $$
    is a reduced crossed-product by an outer action $F$ of $\bbZ$ on $\cO_n^\bbT$ by dualizable bimodules under the group C*-algebra object $\bbC[\bbZ].$
\end{thmalpha}
\noindent Here, $E$ averages the action over $\bbT$ (cf \cite[\S V.4]{MR1402012}), and $\cO_n^\bbT$ is the UHF  algebra $M_{n^\infty} \cong \bigotimes_{m\in\bbN} M_n$. 
We prove this theorem by describing the inclusion's projective quasi-normalizer as the $*$-algebra spanned by the generating isometries $\{s_i\}_{i=1}^n$, and showing explicitly how to recover this inclusion as a crossed-product by a UTC-action. 

\bigskip

This manuscript intends to bridge the gap between the literature in quantum symmetries and  C*-algebras by explaining in detail that many naturally occurring C*-inclusions fit well into the framework of UTCs, while implementing different techniques from those previously utilized.
For instance, our crossed-products by UTC actions do not require passing to the Toeplitz extensions (c.f. \cite{MR2032486}) nor stabilizing \cite[\S2]{MR467330}.

\medskip

To keep this manuscript brief, we follow the notation from \cite{2023arXiv230505072H}, and refer the reader there for the specifics related to C*-correspondences as well as UTCs, and crossed-products by their actions on C*-algebras. 
In {\bf Section \ref{sec:CPAlgebrasDiscrete}} we briefly recall the necessary prerequisites on Cuntz-Pimsner algebras and prove Theorem \ref{thmalpha:CuntzPimsnerCoreDisc}, and in {\bf Section \ref{sec:CuntzCoreDisc}} we develop the proof of Theorem \ref{thmalpha:CuntzCoreDisc}.

\section*{Acknowledgements}
I am grateful to Brent Nelson and Matthew Lorentz who provided tremendous support in the initial stages of this project, and to Pawel Sarkowicz and Jeremy Hume who provided many useful comments. The author was supported by the AMS-Simons Travel Grant 2022 and NSF grants DMS-2001163, DMS-2000331, and DMS-1654159.

\section{Discrete inclusions from Cuntz--Pimsner Algebras}\label{sec:CPAlgebrasDiscrete}

Throughout this section, we let $D$ be a unital C*-algebra, and let $X$ be a $D$-$D$ correspondence, which is assumed to be dualizable (i.e.,  $X\in\fgpBim(D)$),  \emph{faithful} and \emph{full}. Recall faithfullness means that we assume that the left action $D\overset{\rhd}{\hookrightarrow}\End(X_D)$ is given by an injective unital $*$-homomorphism (hence isometric) into adjointable operators.  
Fullness means that $\mathsf{span}\{\langle \eta\ |\ \xi \rangle_{D}\}_{\eta, \xi\in X}$ is dense in $D$. 
In other words, the evaluation map $\ev_{X}:\overline{X}\boxtimes X\to D$ from the duality equations is surjective. 

We start by recording in Lemma \ref{lem:fgpCompacts} that the dualizability/fgp assumption on $X$ forces $\End(X_D)$ and $\End({}_DX)$ to consist of $D$-finite rank operators; that is, every left or right $D$-linear morphism out of $X$ is automatically adjointable and $D$-compact. 
\begin{lem}\label{lem:fgpCompacts}
For $Z\in \fgpBim(D)$, the right $D$-linear endomorphisms (resp. left $D$-linear) satisfy $$\End(Z_D)= \cK(Z_D).$$
\end{lem}
\begin{proof}
Let $u_1,\ldots, u_n\in X$ be a right Pimsner--Popa basis for $Z$. Then for $T\in \End(Z_D)$ and $\xi\in X$ 
    \[
        T(\xi) = T\left( \sum_{i=1}^n u_i \lhd \<u_i \mid \xi \>_D\right)= \sum_{i=1}^n T(u_i) \lhd \<u_i \mid \xi\>_D.
    \]
Hence $T = \sum_{i=1}^n |T(u_i)\rangle\langle u_i | \in \cK(Z_D)$.
\end{proof}

We let
    \[
        \cF(X):= \overline{\bigoplus_{n=0}^\infty X^{\boxtimes_D n}}^{\|\cdot\|_D}
    \]
be the associated Fock space, where by convention $X^{\boxtimes_D 0}:=D$. Then one defines a representation $(\varphi_\infty,\tau_\infty)$ of $X$  on the adjointable right $D$-linear endomorphisms $\End(\cF(X)_D)$ by
    \begin{align*}
        \varphi_\infty(d)(\xi_1\boxtimes \cdots \boxtimes \xi_n) &= (d \rhd \xi_1)\boxtimes\xi_2\boxtimes \cdots \boxtimes \xi_n  &\forall d\in D\\
        \tau_\infty(\xi)(\xi_1\boxtimes \cdots \boxtimes \xi_n) &= \xi\boxtimes \xi_1\boxtimes \cdots \boxtimes \xi_n  &\forall \xi, \xi_1, \cdots, \xi_n\in X.
    \end{align*}
Noting that Lemma~\ref{lem:fgpCompacts} implies the \emph{Katsura ideal} is $J_X=D$ and letting $\pi\colon \End(\cF(X)_D) \to \End(\cF(X)_D)/ \cK(\cF(X)_D)$ be the quotient map, define $\varphi:=\pi\circ \varphi_\infty$ and $\tau:=\pi\circ \tau_\infty$. 
Then $(\varphi,\tau)$ is the universal covariant representation of $X$, and the associated Cuntz--Pimsner algebra for $X$ is defined as 
$$\cO_X:=C^*(\varphi,\tau)\subset \End\left(\cF(X)_D\right)/\cK(\cF(X))$$
(see \cite[Proposition 6.5]{MR2102572}). 
Let $\bbT \overset{\gamma}{\curvearrowright} \cO_X$ be the gauge action: for $d\in D$, $\xi\in D,$ and $z\in \bbT$ we have 
    \[
        \gamma_z(\varphi(d)) = \varphi(d) \qquad \text{ and } \qquad \gamma_z( \tau(\xi)) = z \tau(\xi).
    \]
Letting $\mu$ be the unique Haar measure on $\bbT$ with $\mu(\bbT)=1$,
\begin{align}\label{Dfn:Expectation}
        E_\gamma(x):= \int_{\bbT} \gamma_z(x)\ d\mu(z)
\end{align}
defines a faithful conditional expectation $E_\gamma\colon \cO_X \twoheadrightarrow \cO_X^\gamma$ onto the fixed point subalgebra, which is also known as {\bf the core of} $\cO_X$. We refer to $E_\gamma$ simply by $E$ when no confusion may arise. 

The gauge action of $\bbT$ on $\cO_X$ induces a faithful grading by its spectral subspaces:
\begin{equation}\label{eqn:SpectralGauge}
    \cO_X \cong \overline{\bigoplus_{k\in \bbZ}\cO_X^{\gamma, k}}^{\|\cdot\|},
\end{equation}
diagonalizing the gauge action.

\medskip
We now introduce the main result of this section. 
\begin{thm}[Theorem \ref{thmalpha:CuntzPimsnerCoreDisc}]\label{thm:CuntzPimsnerCoreDisc}
If $D$ is a unital C*-algebra, and $X\in \fgpBim(D)$ is full and faithful, then $$\cO_X^\gamma \overset{E}{\subset} \cO_X$$ is a C*-discrete inclusion.
\end{thm}
\noindent The rest of this section is dedicated to the proof of this result. 

To keep the notation light, we will denote the core $A:=\cO_X^\gamma$, and given $k\in \bbZ$ and $r\in\bbN\cup\{0\}$ with $k+r\geq 0,$ we denote the following  spaces of \emph{right $D$-linear morphisms}: 
\begin{align}\label{eqn:cD-k-r}
    {}^{k+r}\cD^{r}&:= \Hom_{-D}\left(X^{\boxtimes r}\to X^{\boxtimes (k+r)}\right), \text{ and }\\
    {}^{(k)}\cD^\circ &:=\varinjlim\!{}^{\circ}\ {}^{k+r}\cD^{r}.\nonumber
\end{align}
Here, the limit is algebraic with (automatically isometric) inclusion maps of Banach spaces 
\begin{align}\label{eqn:IndLimSys}
    {}^{ k+r}\iota^{r}:{}^{k+r}\cD^{r} &\hookrightarrow {}^{k+r+1}\cD^{r+1}\\
    T&\mapsto T\boxtimes \id_X.\nonumber
\end{align}
(We shall drop the superindices in $\iota$ when no confusion arises.)
For $k\in \bbZ,$ we shall show that (\ref{eqn:IndLimSys}) is canonically an inductive limit system of C*-correspondences (see Lemmas~\ref{lem:Promotion} and \ref{lem:FrobeniusRec}) in the sense of \cite[\S4.1]{MR4717816}, yielding a family of $A$-$A$ C*-correspondences 
$${}^{(k)}\cD:= \varinjlim\ {}^{k+r}\cD^{r}.$$
(This completion procedure exactly matches with the the Banach space inductive limits taken in \cite{MR1010160, MR1638139}, as well as with the operator norm completion of the algebraic limits taken in \cite{MR1658088}.)
If $k=0,$ for each $r$, the space ${}^r\cD^r$ is a unital C*-algebra, and by \cite[Proposition 5.7]{MR2102572} its C*-inductive limit is ${}^{(0)}\cD\cong A=\cO_X^\gamma$.

The proof of Theorem \ref{thm:CuntzPimsnerCoreDisc} will rely on an alternative presentation of $\cO_X$ using the {\bf Doplicher-Roberts picture} \cite{MR1010160} as described in \cite[Proposition 2.5] {MR1658088} and \cite[Proposition 3.2]{MR1638139}:  
\begin{align}\label{eqn:DR-CP}
    \overline{\bigoplus_{k\in \bbZ}\ {}^{(k)}\cD }^{\|\cdot\|}\ &\cong \ \overline{\bigoplus_{k\in \bbZ}\cO_X^{\gamma, k}}^{\|\cdot\|} =\cO_X\\
    X\ni\xi&\mapsto \tau(\xi), \nonumber
\end{align}
as unital $\bbZ$-graded C*-algebras; that is, the isomorphism preserves the spectral subspaces from (\ref{eqn:SpectralGauge}).

In light of (\ref{eqn:DR-CP}), to prove Theorem \ref{thm:CuntzPimsnerCoreDisc} it suffices to show that every term in the direct sum on the left-hand side is finitely generated projective over $A.$ 
In turn, this implies that 
$$\bigoplus_{k\in \bbZ}{}^{(k)}\cD\subseteq \PQN(A\subset \cO_X)$$ 
with their natural unital $*$-algebra structures.

We now consider ${}^{k+r}\cD^{r}$ as a  ${}^{k+r}\cD^{k+r}-{}^{r}\cD^{r}$ bimodule, where the actions are given by post- and pre-composition with the respective endomorphism C*-algebra. 
In fact, with the right inner product 
\begin{equation}\label{eqn:PromotedBimodule}
\langle T\ |\ S\rangle_{{}^{r}\cD^{r}}:= T^* S.    
\end{equation}
these bimodules can be turned into 
right C*-correspondences. 
\begin{lem}\label{lem:Promotion}
    For $k\in\bbZ$ and $r\in\bbN\cup\{0\},$ with $k+r\geq 0$, we have that 
    $$
    {}^{k+r}\cD^{r}\in \rCorr\left({}^{k+r}\cD^{k+r}\ \to \ {}^{r}\cD^{r}\right)\!;
    $$
    that is, ${}^{k+r}\cD^{r}$ is a right ${}^{k+r}\cD^{k+r}-{}^{r}\cD^{r}$ C*-correspondence.
\end{lem}
\begin{proof}
    The inner-product from (\ref{eqn:PromotedBimodule}) is positive definite, since $X\in \rCorr(D\to D)$ is a C*-category, and consequently for every $T\in {}^{k+r}\cD^r$ there exists $R\in {}^r\cD^r$ such that $T^*T = R^*R$.  
    Furthermore, completeness with respect to the induced $\|\cdot\|_{r}$-norm is ensured, since it coincides with the operator norm on ${}^{k+r}\cD^r,$ under which it is already a Banach space. 
\end{proof}

The following is a Frobenius Reciprocity statement which is crucial in our analysis, and we make minor adaptations from \cite[\S2]{{MR1624182}} to make it suitable for us.  
\begin{lem}\label{lem:FrobeniusRec}
    Given $k\in \bbZ$ and $r\in\bbN\cup\{0\}$ with $k+r\geq 0,$ 
    ${}^{k+r}\cD^{r}$ has a canonical structure of a $D$-$D$ bimodule given by 
    $$(d_1\rhd T\lhd d_2) (\xi) = d_1\rhd T(d_2\rhd\xi)\qquad \forall d_1, d_2\in D,\ \forall T\in {}^{k+r}\cD^r,\ \forall \xi\in X^{\boxtimes r}.$$
    Moreover, there is an isomorphism of $D$-$D$ bimodules
    \begin{align*}
        {}_D\left(X^{\boxtimes(k+r)}\boxtimes \overline{X}^{\boxtimes r}\right)_D&\cong\ {}_D({}^{k+r}\cD^{r})_D    \\
        \eta\boxtimes\overline{\xi}&\mapsto |\eta\rangle\langle \xi|. 
    \end{align*}
    Finally, we have 
    $$
    {}^{k+r}\cD^{r}\in \fgpBim({}^{k+r}\cD^{k+r}\to{}^{r}\cD^{r}).
    $$
\end{lem}
\begin{proof}
    The first two claims are exactly \cite[Proposition 2.5]{MR1624182}.
    The last part follows from the $D$-$D$ isomorphism above, since $D\hookrightarrow{}^{r}\cD^{r}$ and $D\hookrightarrow {}^{k+r}\cD^{k+r}$ are the left $D$-actions on $X^{\boxtimes r}$ and $\overline{X}^{\boxtimes (k+r)}$, and $X^{\boxtimes(k+r)}\boxtimes \overline{X}^{\boxtimes r}$ is already $D$-$D$ fgp. 
\end{proof}

The following lemma is a convenient non-degeneracy condition affording us a technical limit tool of \cite{MR4717816}. 
\begin{lem}\label{lem:IL1}
    Given $k\in \bbZ$ and $r\in\bbN\cup\{0\}$ with $k+r\geq 0,$
    $$\iota\left[ {}^{k+r}\cD^r\right]\lhd \left({}^{r+1}\cD^{r+1}\right) =\ {}^{k+r+1} \cD^{r+1}$$
\end{lem}
\begin{proof}
    We use the isomorphism from Lemma \ref{lem:FrobeniusRec}, so that ${}^{k+r}\cD^r \cong X^{k+r}\boxtimes \overline{X}^r.$ 
    Under this map, the inclusion $\iota$ becomes 
    $$
    \iota(\xi_1\boxtimes \cdots \boxtimes \xi_{k+r}\boxtimes \overline{\eta_1}\boxtimes \cdots \boxtimes \overline{\eta_r}) = \sum_{u} \xi_1\boxtimes \cdots \boxtimes \xi_{k+r}\boxtimes u\boxtimes \overline{u}\boxtimes\overline{\eta_1}\boxtimes \cdots \boxtimes \overline{\eta_r},
    $$
    where $\{u\}\subset X$ is a right Pimsner-Popa basis.
    Let
    $$
    \boxtimes_{\ell=1}^{k+r+1}\alpha_\ell\boxtimes \boxtimes_{g=1}^{r+1}\overline{\beta_g} \in X^{k+r+1}\boxtimes \overline{X}^{r+1}$$
    be arbitrary, and we will construct a preimage on $\iota\left[ {}^{k+r}\cD^r\right]\lhd \left({}^{r+1}\cD^{r+1}\right)\!.$
    Since $X$ is full and $D$ is unital, we can choose finitely many vectors $\{x_i, y_i\}_{i\in I}\subset X$ such that $\sum_{i\in I}\langle x_i\ |\ y_i\rangle=1.$
    We now let

    \begin{align*}
    &\sum_{\substack{ i_1, \hdots, i_r \in I}} \iota \left( (\boxtimes_{\ell=1}^{k+r}\alpha_\ell) \boxtimes (\boxtimes_{s=1}^r \overline{x_{i_s}})\right) \lhd \left( (\boxtimes_{t=r}^1 y_{i_t})\boxtimes \alpha_{k+r+1}\boxtimes (\boxtimes_{g=1}^{r+1}\overline{\beta_g}) \right) \\
    &\hspace{-.2cm}=\sum_{\substack{ u, \\ i_1, \hdots, i_r \in I}}  (\boxtimes_{\ell=1}^{k+r}\alpha_\ell) \boxtimes u \boxtimes \overline{u} \boxtimes (\boxtimes_{s=1}^r \overline{x_{i_s}}) \lhd \left( (\boxtimes_{t=r}^1 y_{i_t})\boxtimes \alpha_{k+r+1}\boxtimes (\boxtimes_{g=1}^{r+1}\overline{\beta_g}) \right)\\
    &\hspace{-.2cm}= \sum_{\substack{ u, \\ i_1, \hdots, i_r \in I}} \hspace{-.35cm} (\boxtimes_{\ell=1}^{k+r}\alpha_\ell) \boxtimes\Bigg[ u \lhd \Bigg\langle u \Bigg| \Big\langle x_{i_1} \Big|\ \Big\langle x_{i_2} \Big|\!\cdots\!\langle x_{i_{r-1}} |\ \langle x_{i_r}|\  y_{i_r}\rangle \rhd y_{i_{r-1}}\rangle\!\cdots\!\rhd y_{i_2}\Big\rangle\rhd y_{i_1}\Big\rangle \rhd \alpha_{k+r+1}\Bigg\rangle\boxtimes\\
    &\hspace{14cm} \boxtimes(\boxtimes_{g=1}^{r+1}\overline{\beta_g})\Bigg] \\
    &\hspace{-.2cm}= \sum_{u}  (\boxtimes_{\ell=1}^{k+r}\alpha_\ell)\boxtimes u\lhd \langle u\ |\ \alpha_{k+r+1}\rangle \boxtimes (\boxtimes_{g=1}^{r+1}\overline{\beta_g})  \\
    &\hspace{-.2cm}=  \boxtimes_{\ell=1}^{k+r+1}\alpha_\ell\boxtimes \boxtimes_{g=1}^{r+1}\overline{\beta_g}.
    \end{align*}
\end{proof}

From Lemma \ref{lem:IL1} it then follows that indeed the system of inclusions from (\ref{eqn:IndLimSys}) is an inductive limit system in the precise sense of  \cite[Definition 4.1]{MR4717816}. 
For a fixed $k\in \bbZ$ we take the inductive limit of ${}^{r}\cD^{r}-{}^{k+r}\cD^{k+r}$ C*-correspondences
\begin{equation}\label{eqn:Completion}
{}^{(k)}\cD:=\varinjlim {}^{k+r}\cD^{r}.
\end{equation}
In particular, for $k=0,$ the inductive limit C*-algebra is $A.$ 
Now, by \cite[Lemma 4.3]{MR4717816}, since each ${}^{k+r}\cD^{r}$ is fgp as a ${}^{k+r}\cD^{k+r}-{}^{r}\cD^{r}$ bimodule by Lemma \ref{lem:FrobeniusRec}, then ${}^{(k)}\cD$ is fgp as an $A$-$A$ bimodule, as Pimsner-Popa bases at the finite levels become bases in the limit.
We record this in the following proposition. 
\begin{prop}\label{prop:Fibersfgp}
    For every $k\in \bbZ$ we have 
    $${}^{(k)}\cD \in \fgpBim(A).$$
\end{prop}

\medskip

\begin{proof}[Proof of {\bf Theorem} \ref{thm:CuntzPimsnerCoreDisc}]
    This now follows from Proposition \ref{prop:Fibersfgp}, since 
    $$\bigoplus_{k\in \bbZ}\ {}^{(k)}\cD\ \subseteq\ \PQR\left( A\subseteq \cO_X\right),$$
    and so the projective quasi-normalizer of the inclusion is dense with respect to the C*-norm on $\cO_X.$
\end{proof}

\begin{exs}\label{exs:C,CK}
Note that for $A=\bbC$, one has $\bbC^n \in \fgpBim(\bbC)$ and hence Theorem \ref{thm:CuntzPimsnerCoreDisc} implies the Cuntz algebra $\cO_n:=\cO_{\bbC^n}$ is discrete over its UHF core $\cO_n^\bbT \cong M_{n^\infty}.$
We will further study this inclusion in Section \ref{sec:CuntzCoreDisc}.

More generally, if $D\cong C(\{\sigma\}_{\sigma\in \Sigma})$ is an abelian finite-dimensional C*-algebra with minimal projections $\{p_\sigma\}_{\sigma\in \sigma}$, and $X\in \fgpBim(D)$ is full and faithful, then Theorem \ref{thm:CuntzPimsnerCoreDisc} tells us that Cuntz-Krieger algebras are C*-discrete over their cores. 
Recall, as in \cite[\S1 Example (2)]{MR1426840}, to $X$ corresponds a square matrix $[a_{\sigma, \sigma'}]\in \mathsf{Mat}(\bbN\cup \{0\})$, whose entries are given by $a_{\sigma, \sigma'}=\dim_\bbC( p_{\sigma}\rhd X\lhd p_{\sigma'}).$ 
Here, fullness is equivalent to this matrix having no rows or columns identically zero.
\end{exs}

We now discuss the structure of $A\overset{E}{\subset}\cO_X$ as a crossed-product by a quantum symmetry in lights of Theorem \ref{thm:CuntzPimsnerCoreDisc}, and the characterization of irreducible C*-discrete inclusions from \cite[Theorem A]{2023arXiv230505072H}. 
In general, the inclusion $A\overset{E}{\subset}\cO_X$ will not be irreducible (i.e.,  $A'\cap \cO_X\not\cong \bbC1$), and furthermore, the center $Z(A)= A'\cap A$ might be infinite-dimensional. 
As a consequence, the rigid C*-tensor category $\fgpBim(A)$ is not a UTC, since the endomorphisms of the identity $\End({}_AA_A)\cong Z(A)$ need not be one-dimensional, and so many of the combinatorial or representation-theoretic aspects of UTCs are not readily available. 

Nevertheless, $\cO_X$ can be still expressed as a reduced crossed-product by the action of a unitary tensor category by dualizable $A$-$A$ correspondences under a connected C*-algebra object. 
Namely, the isomorphism from (\ref{eqn:DR-CP}) and Proposition \ref{prop:Fibersfgp} gives a faithful non-full $\dag$-tensor functor 
\begin{align}\label{eqn:FforCP}
    F:\fdHilb(\bbZ)&\to \fgpBim(A)\\
        \bbC_0&\mapsto\ A= {}^{(0)}\cD \nonumber\\
        \bbC_k&\mapsto\ {}^{(k)}\cD \nonumber\\
        \Hom(\bbC_k\to \bbC_\ell)\ni \delta_{k=\ell}\cdot\lambda\id_k&\mapsto (\delta_{k=\ell}\cdot\lambda \id_{{}^{(k)}\cD}) \in \Hom_{A-A}\left({}^{(k)}\cD\to {}^{(\ell)}\cD\right)\!.\nonumber
\end{align}
With unitary tensorator defined by 
\begin{align}
    F^2_{k,\ell}:{}^{(k)}\cD\boxtimes_A {}^{(\ell)}\cD&\to {}^{(k+\ell)}\cD\\
    T\boxtimes S&\mapsto T\circ S.\nonumber
\end{align}
Here $T$ and $S$ are assumed elements in ${}^{(k)}\cD^\circ$ and ${}^{(\ell)}\cD^\circ$, and by tensoring with sufficient identities on the right---an isometric transformation---we make sense of the composite by embedding them at the same level. 
This map is well-defined, as it is directly seen to factor through the $A$-balanced tensor product.
Since at algebraic levels, all $D$-linear maps are compact, it follows that $F^2_{k,\ell}$ has dense range.
Unitarity then follows from an immediate diagrammatic computation with the $A$-valued inner product. Therefore the map extends to a unitary in the limit for each $k, \ell$. 

Compare the action in (\ref{eqn:FforCP}) with that of (\ref{eqn:F}) associated to the Cuntz algebra $\cO_n$, where outerness (i.e.,  functor is full and faithful) is afforded by the orthogonality of the generators, which does not hold in larger generality. 
The action from (\ref{eqn:FforCP}) together with the connected C*-algebra object $\oplus_{k\in \bbZ}\bbC_k=\bbC[\bbZ]\in \Vec(\fdHilb(\bbZ))$ gives: 
$$\bigoplus_{k\in \bbZ}F(k)\otimes \bbC[\bbZ](k)= \bigoplus_{k\in \bbZ}\ {}^{(k)}\cD,$$
with its structure as a unital $*$-algebra. 
And thus, $\cO_X\cong A\rtimes_{r,F}\bbC[\bbZ],$
in the sense of \cite[\S3]{2023arXiv230505072H}, by the isomorphism from Equation \ref{eqn:DR-CP}. 
We record this result in the following corollary. 
\begin{cor}\label{cor:CPCrossedProduct}
    Let $D$ be a unital C*-algebra, $X\in \fgpBim(A)$ be full and faithful, and $F$ is the unitary tensor functor from (\ref{eqn:FforCP}).
    We then have an expectation-preserving $*$-isomorphism of C*-discrete inclusions
           $$\left(A\overset{E}{\subset} \cO_X\right)\cong \left(A\overset{E'}{\subset}A\rtimes_{r,F}\bbC[\bbZ]\right)\!.$$ 
    Here, $E= E_\gamma$ is the canonical expectation onto the fixed points of the gauge action, and $E'$ is the canonical expectation projecting onto the $0$-graded component.
\end{cor}

\section{Cores of Cuntz Algebras are C*-Discrete}\label{sec:CuntzCoreDisc}
A well-known inclusion of simple C*-algebras is given by taking the fixed point subalgebra of the \emph{gauge action} on the Cuntz algebra in $n\in \bbN$ generators, yielding the subalgebra $n^\infty$-UHF C*-algebra $M_{n^\infty}$ \cite[\S V.4]{MR1402012} \cite{Rae05}. 
Namely, for $n\geq 2$, $\cO_n$ is the universal C*-algebra generated by  $\{s_i\}_{i=1}^n$ with relations $\sum_{i=1}^ns_is_i^* =1$ and for each $i$ and $j$ one has $s_i^*s_j = \delta_{i=j}1.$ 
Equipped with the universal property of this construction, one can construct a group homomorphisms $\bbT\to \Aut(\cO_n)$ given by $t\cdot s_j= e^{2\pi i t}s_j$ for $t\in [0,1].$ 
The fixed-point C*-algebra $A = \cO_n^\bbT$ is the span of words of the form $s_\nu s_\mu^*$, where $\nu = (\nu_1, \nu_2,...,\nu_{|\nu|})$ and $\mu = (\mu_1,\mu_2,..., \mu_{|\mu|})$ are finite sequences in $\{1,2,\hdots, n\}$ with lengths $|\nu|=|\mu|$, which in turn yields $A\cong M_{n^\infty}.$ 
Our goal in this section is to give a detailed and direct description of the inclusion $M_{n^\infty}\overset{E}{\subset}\cO_n$, where $E$ is the faithful conditional expectation determined by
$$
E(s_\nu s_\mu^*) = \delta_{|\nu|=|\mu|} s_\nu s_\mu^*,
$$
from which we shall conclude the inclusion is irreducible C*-discrete:
\begin{thm}[Theorem \ref{thmalpha:CuntzCoreDisc}]\label{thm:CuntzCoreDisc}
For any $n\geq 2,$ the inclusion $M_{n^\infty}\overset{E}{\subset}\cO_n$ arising from the fixed points of the canonical gauge action is irreducible and C*-discrete.
\end{thm}

Many facts about the inclusion $M_{n^\infty}\subset \cO_n$ are well-established, like the fact that it is \emph{C*-irreducible} \cite{MR4599249, 2021arXiv210511899Raddendum}; that is, every intermediate unital subalgebra $M_{n^\infty}\subseteq D\subseteq \cO_n$ is simple, and that $\cO_n = C^*(M_{n^\infty}, s_1)=\cO_n\cong M_{n^\infty}\rtimes_{\lambda_1} \bbN$, where the latter means that $\cO_n$ can be obtained as the ``crossed-product'' of $A$ by the semi-group of non-negative integers by the (non-unital) $*$-endomorphism $\lambda_1 = \Ad(s_1)\in\End(\cO_n)$, allowing us to view elements in $\cO_n$ as \emph{Fourier series with coefficients in $A$}. 
We shall make this crossed-product structure explicit and express it in our terms in Corollary  \ref{cor:OnCrossedProduct}. 
To ease notation, in the following we write $A$ instead of $M_{n^\infty}$. 

\medskip
We shall first show that $\PQN(A\subset \cO_n)$ contains the span of all finite words $s_\nu s_\mu^*,$ so that $A\subset \cO_n\in\CDisc.$
\begin{construction}
Consider the following unital $*$-preserving maps on $A:$
\begin{align*}
    &\lambda:A\hookrightarrow\lambda[A]\subsetneq A\qquad \qquad &&\cE:A\twoheadrightarrow \cE[A] &&& \cL:A\twoheadrightarrow \cL[A]=A\\
    &\hspace{.7cm}a\mapsto \sum_{i=1}^n s_i as_i^*, &&\hspace{.7cm} a\mapsto \sum_{i=1}^n s_is_i^* as_i s_i^*, &&& a \mapsto \sum_{i=1}^ns_i^*as_i.
\end{align*}
One can easily check that $\lambda$ is a unital $*$-endomorphism and is thus injective as $A$ is simple, that $\cE$ is $\lambda[A]$-$\lambda[A]$ bimodular conditional expectation; i.e., , faithful ucp, that $\cL(s_\nu s_\mu^*)=\delta_{\nu_1 = \mu_1}s_{\nu_2,\hdots\nu_{|\nu|}}s_{\mu_2,\hdots\mu_{|\nu|}}^*$, that $\cL\circ\lambda = n\id_A$ and that for all $a, a'\in A$ we have $\cL(\lambda(a)a') = a\cL(a').$  
Notice, $\cL$ is a \emph{transfer operator} as described in \cite[2.1 Definition]{MR2032486} but our $\cE$ is slightly different from his $E(a):=\sum_{i,j=1}^ns_is_jas_j^*s_i$, yielding bimodules whose structures only differ by the scalar $n1$. By \cite[2.3 Proposition]{MR2032486} $E[A]=\lambda[A],$ so $\cL$ is \emph{non-degenerate}. 

Automatically, we now have that ${_\lambda}A_A\in \rCorr(A\to A)$ is a right C*-correspondence with the trivial right Hilbert $A$-module structure, and for $a\in A,$ the left $A$ action $a\rhd$ is given by left multiplication by $\lambda(a)$ for $a\in A$.  
In the following, using the conditional expectation $\cE$ and the \emph{transfer operator} $\cL$ we show the much more stronger statement that in fact ${_\lambda}A_A\in\fgpBim(A).$

We now impart ${_\lambda}A$ with the $A$-valued form given by: 
$${_A}\langle \eta, \xi\rangle:= \cL\circ\cE(\eta\xi^*).$$
For each $\eta,\xi\in{_\lambda}A$ and each $a\in A$ we have that ${_A}\langle \eta, \xi\rangle^* = (\cL\circ\cE(\eta\xi^*))^* = \cL\circ\cE(\xi\eta^*) =  {_A}\langle \xi, \eta\rangle,$ so the form is anti-symmetric. 
From the identity $\cL(\lambda(a)\cE(b)) = a\cL(b) = a\cL(\cE(b))$ for each $a,b\in A$ we obtain that ${_A}\langle a\rhd\eta, \xi\rangle = \cL\circ\cE(\lambda(a)\eta\xi^*) = a\cL(\cE(\eta\xi^*)) = a{_A}\langle \eta, \xi\rangle$, so that the form is left $A$-linear. 
Furthermore, it can be seen directly from the definition of the maps that ${_A}\langle \eta, \eta\rangle\geq 0 $ and ${_A}\langle \eta, \eta\rangle = 0$ if and only if $\eta=0.$ 
Therefore $({_\lambda}A, {_A}\langle \cdot, -\rangle)$ is a left Hilbert $A$-module. 
By noticing that for each $a\in A$ we have ${_A}\langle \eta, \xi\lhd a\rangle = {_A}\langle \eta, \xi a\rangle = \cL\circ\cE(\eta a^*\xi^*) =  {_A}\langle \eta\lhd a^*, \xi\rangle$ we conclude that the right $A$-action on ${_\lambda}A$ is by left $A$-adjointable operators. 
Consequently, $({_\lambda}A_A, {_A}\langle \cdot, -\rangle)$ is also a left $A$-$A$ correspondence. Finally, for each $a\in A$
    \[
       \|a\| = \left\|\sum_{i=1}^n s_is_i^*a\right\|  \leq  \sum_{i=1}^n \| s_i^* a\| \leq \sum_{i=1}^n \left\| \sum_{j=1}^n s_j^* aa^* s_j \right\|^{\frac12} \leq \sum_{i=1}^n \left( \sum_{j=1}^n \|a\|^2 \right)^{\frac12},
    \]
which implies $\|\xi\|_A \leq n\,  {}_A \|\xi\| \leq n^{\frac32} \|\xi\|_A$ for all $\xi\in {}_\lambda A_A$. Thus the topologies induced by the left and right $A$-valued inner products agree and  therefore ${}_\lambda A_A$ is a Hilbert $A$-$A$ bimodule. 
Notice the middle inequality above follows from the positivity condition $s_i^*aa^*s_i\leq \sum_{j=1}^n s_j^*aa^*s_j$ for each $i$.

\begin{lem}\label{lem:lambdaA}
    The Hilbert $A$-$A$ bimodule $({_\lambda}A_A, {_A}\langle \cdot, -\rangle, \langle-\mid\cdot\rangle_A)\in \fgpBim(A)$ with finite left $A$ Pimsner--Popa basis given by $\{s_is_j^*\}_{i,j=1}^n$ and finite right $A$ Pimsner--Popa basis $\{1\}.$ 
    Furthermore, ${_\lambda}A_A$ is balanced (i.e.,  $\ell-\mathsf{Ind}_W({_\lambda}A) = n1 = r-\mathsf{Ind}_W(A_A)$) with Watatani index $\mathsf{Ind}_W({_\lambda}A_A)=n^21.$
\end{lem}
\begin{proof}
    It only remains to show the statement about the left $A$ basis. 
    For an arbitrary word of the form $s_\nu s_\mu^*$ with $|\nu| = |\mu|\in\bbN$ we have that
    \begin{align*}
        \sum_{i,j=1}^n {_A}\langle s_\nu s_\mu^*, s_js_i^*\rangle\rhd s_js_i^*       &= \sum_{i,j=1}^n \lambda\circ\cL\circ\cE(s_\nu s_\mu^* s_is_j^*)s_js_i^*\\
        &= \sum_{i,j=1}^n \delta_{i=\mu_1}\lambda\circ\cL\circ\cE(s_\nu s_{\mu_2\hdots\mu_{|\mu|}}^* s_j^* )s_js_i^*\\
        &= \sum_{j=1}^n \lambda\circ\cL\circ\cE(s_\nu s_{\mu_2\hdots\mu_{|\mu|}}^* s_j^*)s_j s_{\mu_1}^*\\
        &= \sum_{j,k=1}^n \lambda\circ\cL (s_k\underbrace{s_k^* s_\nu}_{\delta_{k = \nu_1}} s_{\mu_2\hdots \mu_{|\mu|}}  \underbrace{s_j^*s_k}_{\delta_{k=j}}s_k^*)s_j s_{\mu_1}^*\\ 
        &= \lambda\circ\cL(s_\nu s_{\mu_2\hdots\mu_{|\mu|}}^* s_{\nu_1}^*)s_{\nu_1}s_{\mu_1}^*\\
        & = \lambda(s_{\nu_2\hdots\nu_{|\nu|}} s_{\mu_2\hdots\mu_{|\mu|}}^*) s_{\nu_1}s_{\mu_1}^*\\
        &= \sum_{\ell=1}^n s_\ell s_{\nu_2\hdots\nu_{|\nu|}} s_{\mu_2\hdots\mu_{|\mu|}}^* s_\ell^* s_{\nu_1}s_{\mu_1}^* = s_\nu s_\mu^*. 
    \end{align*}
    Thus any $a\in A$ which is a finite linear combination of words as above can be expanded using the left $A$-basis, and since the span of words is dense in $A$ this establishes that ${_\lambda}A_A\in\fgpBim(A)$. 
    
    The Watatani index is computed using \cite[Equation 1]{2023arXiv230505072H} so we have that $\mathsf{Ind}_W({_\lambda}A_A)= r-\mathsf{Ind}(A_A)\cdot \ell-\mathsf{Ind}({_\lambda}A) = {_A}\langle 1, 1\rangle \cdot \sum_{i,j=1}^n \langle s_js_i^*\mid s_js_i^* \rangle{_A} =  n1\cdot \sum_{i,j=1}^n  s_is_j^* s_js_i^* = n^2 1.$
\end{proof}

\begin{prop}\label{prop:lambdaA}
    With the conventions from from Lemma \ref{lem:lambdaA} we have the following: 
    \begin{enumerate}[label=(\arabic*)]
        \item The left Pimsner--Popa basis $\{s_js_i^*\}_{i,j=1}^n$  is orthonormal; i.e.,  ${_A}\langle s_js_i^*, s_ks_\ell^* \rangle= \delta_{j=k}\delta_{i=\ell}1,$ and we can express 
        ${_\lambda}A\cong \bigoplus_{i,j=1}^n (A\rhd s_js_i^*)$ as left Hilbert $A$-modules, where each block $(A\rhd s_js_i^*)\in \fgpMod(A)$ has left-Watatani index equal to $s_is_i^*.$
        \item For each $1\leq i\leq n,$ we have an irreducible fgp bimodule  ${_\lambda}((s_is_i^*)[A])_A\subset {_\lambda}A_A$ yielding the orthogonal decomposition ${_\lambda}A_A \cong \bigoplus_{i=1}^n{_\lambda}((s_is_i^*)[A])_A,$ with $\ell-\mathsf{Ind}_W(s_is_i^*A) = 1 = r-\mathsf{Ind}_W(s_is_i^*A).$  
        \item We have that the $A$-$A$ bimodular (automatically adjointable) endomorphisms of ${}_\lambda A_A$ satisfy: $\End^\dag({_\lambda}A_A) \cong \lambda[A]'\cap A \cong M_n(\bbC).$
        \item For each $i=1,\hdots n,$ any non-zero $\xi\in (s_is_i^*)A$ algebraically generates the entire bimodule; 
        that is, with the notation from \cite[Equation (5)]{2024arXiv240918161H}  we have  $((s_is_i^*)[A])^{\diamondsuit\diamondsuit} = (s_is_i^*)[A]$. 
    \end{enumerate}
\end{prop}
\begin{proof}
(1) follows from direct computation.

(2): we have for each $i$ the commutation relation $\lambda(a)s_i = s_i a$ for every $a\in A$ in combination with the simplicity of $A$ yields 
$$\overline{A\rhd s_is_i^*\lhd A}^{\|\cdot\|_A}=\overline{\lambda[A] s_is_i^*A}^{\|\cdot\|_A}= \overline{(s_is_i^*)(\lambda[A]A)}^{\|\cdot\|_A} = (s_is_i^*)[A] \in \fgpBim(A)$$
with finite left-$A$ Pimsner--Popa basis $\{s_is_j^*\}_{j=1}^n.$  
We shall now establish the irreducibility of $(s_is_i^*)[A]$ by computing its left and right Watatani indices. 
We have that 
\begin{align*}
    \ell-\mathsf{Ind}({}_\lambda s_is_i^*A) &= \sum_{j=1}^n\langle s_is_j^*\ |\ s_i s_j^* \rangle_A =1,\ \text{ and } \\
    r-\mathsf{Ind}(s_is_i^*A_A) &=\ {}_A\langle s_is_i^*, s_is_i^*\rangle= \sum _{j,k=1}^n s_k^*  s_js_j^*    s_is_i^*         s_js_j^*   s_k=1.
\end{align*}
And thus,  $\Ind_W({}_\lambda s_is_i^*A_A)=1$ is irreducible. The orthogonality of $s_is_i^*A$ and $s_js_j^*A$ with respect to left and right inner products follows from a quick computation. The direct sum decomposition follows at once from $\sum_1^n s_is_i^* =1_A.$

(3):  
Notice first that the matrix units  $\{s_js_i^*\}_{i,j=1}^n\subset\lambda[A]'\cap A.$
In fact, for each $i,j=1,\cdots, n$, the map 
\begin{align*}
    M_{s_js_i^*}:\ {}_\lambda(s_is_i^*)A_A&\to\ {}_\lambda(s_js_j^*)A_A\\
    \xi&\mapsto s_js_i^*\xi
\end{align*}
is readily seen to be an $A$-$A$ bilinear unitary map, and so 
$$
\End^\dag({}_\lambda A_A) \cong \lambda[A]'\cap A \cong \mathsf{span}_\bbC\{s_js_i^*\}_{i,j=1}^n\cong M_n(\bbC).
$$

(4): Let $\xi=s_is_i^* a\in (s_is_i^* )A$ be non-zero with $a\in A$. For each $j=1,\ldots, n$ define $b_j:=s_i^* as_j\in A$ so that
    \[
        \xi = \sum_{j=1}^n s_is_i^* a s_js_j^* = \sum_{j=1}^n s_i b_j s_j^*.
    \]
Let $1\leq j\leq n$ be such that $\xi s_j s_j^* = s_i b_j s_j^*$ is non-zero, which exists lest $\xi=0$. Since
    \[
        \lambda[A] \xi A \supseteq \lambda[A] \xi s_js_j^* A = \lambda[A] s_i b_j s_j^* A,    
    \]
it suffices to show $s_ib_js_j^*$ algebraically generates $(s_is_i^*)A$. As $A$ is unital and simple, we have $Ab_jA=A$ and thus
    \[
        \lambda[A] s_i b_j s_j^* A = \lambda[A] s_i b_j s_j^* \lambda[A]A = s_i A b_j A s_j^* A = s_i A s_j^* A = s_i s_j^* \lambda[A] A = s_i s_i^* A.\qedhere
    \]
\end{proof}
\end{construction}

In lights of Proposition \ref{prop:lambdaA}, it is meaningful to consider the irreducible subbimodules of  ${_\lambda}A_A$, whose structure is easily seen to be described as follows: For each $i=1,...,n$ we have the (non-unital) $*$-preserving maps:
\begin{align*}
    &\lambda_i:A\hookrightarrow\lambda_i[A]\subsetneq (s_is_i^*)A\qquad \qquad &&\cE_i:A\twoheadrightarrow \cE_i[A] &&& \cL_i:A\twoheadrightarrow \cL_i[A]\\
    &a\mapsto  s_i as_i^*, && a\mapsto  s_is_i^* as_i s_i^*, &&& a \mapsto s_i^*as_i.
\end{align*}
So, for each $i=1,...,n$ we obtain fgp bimodules ${_{\lambda_i}}(s_is_i^*)[A]_A$ with the obvious structure as right Hilbert $A$-modules, left action of $a\in A$ by left multiplication by $\lambda_i(a),$ left $A$-valued inner product given by 
$$
{_A}\langle \eta,\xi\rangle^i:= \cL_i\circ\cE_i(\eta\xi^*)=\cL_i(\eta\xi^*).
$$ 
Evidently, the left $A$ action through $\lambda$ on each corner ${_{\lambda}}s_is_i^*[A]_A$ restricts to $\lambda_i,$ and similar statements hold for the $A$-valued inner products, and thus ${_\lambda}A_A \cong \bigoplus_1^n {_{\lambda_i}}(s_is_i^*)A_A.$ 
Left multiplication by $s_is_1^*$ is a bimodule isomorphism  ${}_{\lambda_1}(s_1s_1^*)A_A\cong {}_{\lambda_i}(s_is_i^*)A_A$ for all $i=1,\hdots, n.$

Notice the relative tensor product yields the composition bimodule  ${_{\lambda_1\circ\lambda_1}}(s_1s_1s_1^*s_1^*)[A]_A$ with trivial right $A$ structure, but on the left the $A$-action is twisted by the $*$-homomorphism $\lambda_1\circ\lambda_1,$ with left $A$-valued product ${_A}\langle\eta,\xi\rangle^{11}:= \cL_{11}(\eta\xi^*):=(s_1s_1)^*(\eta\xi^*)(s_1s_1).$ 
We now record how these pieces interact: 

\begin{prop}\label{prop:lambdaFusRules}
    For each $m\in\bbN$ we have the following $A$-$A$ unitary maps: 
    \begin{align*}
        L_{(s_1^*)^m}: {_{\lambda_1^m}}s_1^m(s_1^*)^mA &\cong\ {}_A(s_1^*)^mA_A\\
        \xi&\mapsto (s_1^*)^m\xi, \\
        \text{ as well as }\hspace{4cm} &\ \\
        R_{s_1^m}: \overline{{_{\lambda_1^m}}s_1^m(s_1^*)^mA} &\cong\ {}_A{As_1^m}_A\\
        \overline{\xi}&\mapsto \xi^*s_1^m.
    \end{align*}
\end{prop} 
\begin{proof}
Straightforward computation. 
\end{proof}

We now bring the $A$-$A$ bimodules from Section \ref{sec:CPAlgebrasDiscrete} into this framework: 
\begin{lem}\label{lem:GradedBimods}
    For each $k\in \bbZ$ we have a unitary isomorphism of $A$-$A$ bimodules
    $$
    {}^{(k)}\cD \cong \begin{cases}
        As_1^k & \text{if } k >0,\\
        A & \text{if } k=0,\\
        (s_1^*)^{-k}A & \text{if } k <0.\\
    \end{cases}
    $$
\end{lem}
\begin{proof}
The case $k=0$ was established above Equation \ref{eqn:DR-CP}. 
Fix $k\in \bbN.$ 
Taking $s_vs_\mu^*\in A$ with words $\nu$ and $\mu$ of length $m$ we see that 
$$
s_\nu s_\mu^*s_1^k = 
\begin{cases}
\delta_{\mu_{m-k+1}=\cdots =\mu_m =1}\cdot s_\nu s_{\mu_1, \hdots, \mu_{m-k}}^* & \text{ if } k\leq m,\\
\delta_{\mu_{1}=\cdots =\mu_m =1}\cdot s_\nu s_1^{k-m} & \text{ if } k> m.
\end{cases}
$$
In either case, the isomorphism $\tau$ from Equation \ref{eqn:DR-CP} when restricted to the $k$-graded component maps the elementary operators 
\begin{align*}
{}^{(k)}\cD^\circ\ni |u_{\nu_1}\boxtimes\cdots \boxtimes u_{\nu_{r+k}}\rangle \langle u_{\mu_1}\boxtimes \cdots \boxtimes u_{\mu_r}|\ &\overset{\tau}{\mapsto}\ s_{\nu_1}\hdots s_{\nu_{r+k}}s_{\mu_1}^*\hdots s_{\mu_{r}}^*\in \cO_n^{\gamma, k},\\
{}^{(k)}\cD^\circ\ni u_{\nu_1}\boxtimes\cdots \boxtimes u_{\nu_{r-k}} \boxtimes  u_{{1}}^{\boxtimes k} &\overset{\tau}{\mapsto}\ s_{\nu_1}\hdots s_{\nu_{r-k}}s_{1}^k \in \cO_n^{\gamma, k}, 
\end{align*}
where $\{u_\ell\}$ denotes a right Pimsner-Popa basis for $X\cong \bbC^n$; i.e.,  an orthonormal basis. 
Consequently, ${}^{(k)}\cD^\circ$ (recall Equation \ref{eqn:cD-k-r}) is isomorphic to the algebraic linear span of words of the forms above. 
Since $As_1^k$ is densely linearly spanned by words of the same forms, the conclusion follows applying $\tau$ to the closed linear span of elementary tensors in ${}^{(k)}\cD^\circ$. 

The case $k<0$ is similar and unitarity is manifest from definitions.
\end{proof}

\medskip

We denote by $\cB$ the $\cO_n$-$A$ right C*-correspondence given by the completion of $B:=\cO_n$ in the right $A$-norm $\|b\|_A^2 = \|E(b^*b)\|$ for $b\in\cO_n.$ That is, $\cB = \overline{\cO_n\Omega}^{\|\cdot\|_A}.$

\begin{prop}\label{prop:GeneratorsInPQN}
    We have that 
    $$
    \{s_i^*\}_{i=1}^n \subset \PQN(A\subset \cO_n)
    $$
    and consequently
    $$
        *-\mathsf{Alg}\{s_i\}_1^n \subseteq \PQN(A\subset \cO_n). 
    $$
\end{prop}
\begin{proof}
    We prove the first statement for $s_1^*$ as the other are analogous. 
    It suffices to see that the $A$-$A$ bimodule generated by $s_1^*$ is fgp, as it clearly lies within $\cO_n$ as $A$ is a subalgebra. 
    We have 
    $$
        A\rhd s_1^* \lhd A = A s_1^*A = s_1^*\lambda[A]A = s_1^*A\in \fgpBim(A),
    $$
    since $\lambda$ is unital and we have used Lemma \ref{lem:GradedBimods} to see that the algebraic bimodule we generated is already closed. This proves the first claim. 

    The second follows at once from \cite[Lemma 4.12]{2023arXiv230505072H}, as the projective quasi-normalizer is a $*$-subalgebra intermediate to $A$ and $\cO_n. $
\end{proof}

\begin{remark}\label{rmk:CstarIrred}
    In \cite{2021arXiv210511899Raddendum}, R{\o}rdam shows that $A\subset \cO_n$ is \emph{C*-irreducible} (i.e.,  all intermediate subalgebras are simple) using that $E$ has the \emph{pinching property} \cite[Definition 3.13]{MR4599249}. 
\end{remark}

We are now in position to prove Theorem \ref{thm:CuntzCoreDisc}:
\begin{proof}{\bf (Theorem \ref{thm:CuntzCoreDisc})}
    That $A\subset \cO_n$ is C*-discrete follows from Proposition \ref{prop:GeneratorsInPQN}, since $\{s_i\}_{i=1}^n$ generate $\cO_n$ as a C*-algebra, and so $\PQN(A\subset \cO_n)$ is dense in norm in $\cO_n.$

    Irreducibility of this inclusion is granted by Remark \ref{rmk:CstarIrred}.
\end{proof}

\medskip 

By \cite[Theorem 4.3]{2023arXiv230505072H} we know that $A\subset \cO_n$ is a reduced crossed-product inclusion $A\subset A\rtimes_{F,r} \bbB.$ 
Here, $F:\cC_{A\subset\cO_n}\to \fgpBim(A)$ is an outer UTC action of the support category $\cC_{A\subset\cO_n}$ generated by the $A$-$A$ bimodules of the inclusion, and $\bbB\in\rCorr\left(\cC_{A\subset\cO_n}^{\op}\right)$ is the connected C*-algebra object given by 
\begin{align}\label{eqn:BAlgObj}
\bbB(K):= \rCorr_{A-A}(K\to \cB)^\diamondsuit = \rCorr_{A-A}(K\to \cB)
\end{align}
(c.f. \cite[Theorem 2.7]{2023arXiv230505072H}).  
We shall now describe $\cC_{A\subset O_n},$ $F$ and $\bbB$ in full detail.

\begin{prop}\label{prop:PropertiesDk}
For $k\in \bbZ$, the $A$-$A$ bimodules ${}^{(k)}\cD$ are irreducible and pairwise non-isomorphic, with dual 
$$\overline{{}^{(k)}\cD}\cong {}^{(-k)}\cD.$$ 
Moreover, for $k,\ell\in\bbZ$ we have
$$
{}_A{}^{(k)}\cD\boxtimes_A {}^{(\ell)}\cD_A \cong\ {}_A{}^{(k+\ell)}\cD_A.
$$
\end{prop}
\begin{proof}
    Let $k\in\bbZ$ be arbitrary.
    We first show that $\Hom_{A-A}(A\to {}^{(k)}\cD)\cong \delta_{k=0}\cdot\bbC1$.
    By Lemma \ref{lem:GradedBimods}, we have that ${}^{(k)}\cD = As_1^k,$ assuming $k>0.$ 
    If $T:A\to {}^{(k)}\cD$ is nonzero, then there is $t\in A$ with $aT(1)=T(a)=ts_1^ka=T(1)a,$ for all $a\in A.$
    Therefore, $ts_1^k\in \cO_n\cap A'= \bbC1$. 
    Without loss of generality, we assume $ts_1^k=1\in A.$
    Therefore, we have 
    $$\delta_{k=0}1=E(ts_i^k)=E(1)=1,$$
    where $E$ is the canonical faithful expectation from (\ref{Dfn:Expectation}) that projects onto the $0$-graded component. 
    Necessarily, $k=0$, lest $T=0.$ The case where $k\leq 0$ is dealt with similarly. 

    The map $\overline{As_1^k}\to (s_1^*)^kA$ given by $\xi s_1^k\mapsto (\xi s_1^k)^*$ is readily seen to be an $A$-$A$ unitary isomorphism, proving $\overline{{}^{(k)}\cD} \cong {}^{(-k)}\cD$. 

    We now describe the fusion rules for the ${}^{(k)}\cD$ bimodules. 
    For $k\in\bbZ$, ${}^{(k)}\cD$ is the norm-closure of the linear span of the words
    described in the proof of Lemma \ref{lem:GradedBimods}.
    Considering the multiplication map 
    \begin{align*}
        {}^{(k)}\cD\boxtimes_A {}^{(\ell)}\cD &\to \cO_n\\    
        \eta\boxtimes\xi&\mapsto \eta\xi,
    \end{align*}
    we see at once this is a well-defined $A$-$A$ bimodular map whose range is  ${}^{(k+\ell)}\cD.$

    For $k,\ell\in\bbZ$, let $T\in\Hom_{A-A}({}^{(k)}\cD\to {}^{(\ell)}\cD)$. Then, by Frobenius Reciprocity, there corresponds a map $T'\in\Hom_{A-A}(A\to {}^{(\ell-k)}\cD).$ 
    As discussed above, $T'\neq 0$ if and only if $k=\ell.$ 
    In this case, then $T\in \bbC\id$ as $T'\in A\cap A'=\bbC 1.$
    Therefore, the ${}^{(k)}\cD$ are pairwise non-isomorphic irreducible $A$-$A$ bimodules.
\end{proof}

\begin{cor}\label{cor:OnCrossedProduct}
    For each $n\geq 2$ we have an expectation preserving $*$-isomorphism of irreducible C*-discrete inclusions:
        $$\left(A\overset{E}{\subset} \cO_n\right)\cong \left(A\overset{E'}{\subset}A\rtimes_{F,r}\bbC[\bbZ]\right).$$ 
    Here, $E$ is the canonical expectation onto the fixed points of the gauge action, and $E'$ is the canonical expectation that maps formal sums with coefficients in $A$ to the component graded by $0\in\bbZ.$ The unitary fully-faithful tensor functor $F$ is described in (\ref{eqn:F}). 
\end{cor}
\begin{proof}
From Proposition \ref{prop:PropertiesDk} we conclude that the category of $A$-$A$ bimodules supported on the inclusion $A\subset \cO_n$ is 
$$
\cC_{A\subset \cO)n} = \left\langle\  \left\{{}^{(k)}\cD   \right\}_{k\in\bbZ},\ \oplus,\ \boxtimes_A,\ \overline{\ \cdot\ }\right \rangle \cong \fdHilb(\bbZ),
$$

By Proposition \ref{prop:PropertiesDk} there is an outer action (i.e.,  fully-faithful unitary tensor functor) 
\begin{align}\label{eqn:F}
\begin{split}
    F:\fdHilb(\bbZ)&\to \cC_{A\subset\cO_n}\subset\fgpBim(A)\\
    \bbC_0 &\mapsto {_A}A_A,\\
    \bbC_k&\mapsto\ {}^{(k)}\cD \cong As_1^k,\\ 
    \bbC_{-k}&\mapsto \ {}^{(-k)}\cD \cong (s_1^*)^kA,\hspace{2cm} \text{ for } k\in \bbN, 
\end{split}
\end{align}

The C*-algebra object from Equation (\ref{eqn:BAlgObj}) is the group algebra object $\bbC[\bbZ]$ internal to $\fdHilb(\bbZ)$ given by 
$$
\bbC[\bbZ](\ell):=\Hom\left(\bbC_\ell \to \bigoplus_{k\in\bbZ}\bbC_k\right) \cong \bbC. 
$$

From Equation \ref{eqn:DR-CP} we explicitly recover the Cuntz algebra as a reduced crossed-product by an outer action of a tensor category under a connected C*-algebra object: 
$$
    A\rtimes_{r,F}\bbC[\bbZ] = \overline{\bigoplus_{k\in \bbZ}} F(k)\otimes \bbC[\bbZ](k)
    = \overline{\bigoplus_{k\in \bbZ}}\ {}^{(k)}\cD\otimes \bbC\cong \cO_n.
$$
\end{proof}

\bibliographystyle{amsalpha}
\bibliography{bibliography}

\providecommand{\bysame}{\leavevmode\hbox to3em{\hrulefill}\thinspace}
\providecommand{\MR}{\relax\ifhmode\unskip\space\fi MR }
\providecommand{\MRhref}[2]{%
  \href{http://www.ams.org/mathscinet-getitem?mr=#1}{#2}
}
\providecommand{\href}[2]{#2}
\begin{thebibliography}{CHPJ24}

\bibitem[AEE98]{MR1467459}
Beatriz Abadie, S{\o}ren Eilers, and Ruy Exel, \emph{Morita equivalence for
  crossed products by {H}ilbert {$C^*$}-bimodules}, Trans. Amer. Math. Soc.
  \textbf{350} (1998), no.~8, 3043--3054. \MR{1467459}

\bibitem[CHPJ24]{MR4717816}
Quan Chen, Roberto Hern\'{a}ndez~Palomares, and Corey Jones, \emph{K-theoretic
  classification of inductive limit actions of fusion categories on
  {AF}-algebras}, Comm. Math. Phys. \textbf{405} (2024), no.~3, Paper No. 83,
  52. \MR{4717816}

\bibitem[Cun77]{MR467330}
Joachim Cuntz, \emph{Simple {$C\sp*$}-algebras generated by isometries}, Comm.
  Math. Phys. \textbf{57} (1977), no.~2, 173--185. \MR{467330}

\bibitem[Dav96]{MR1402012}
Kenneth~R. Davidson, \emph{{$C^*$}-algebras by example}, Fields Institute
  Monographs, vol.~6, American Mathematical Society, Providence, RI, 1996.
  \MR{1402012}

\bibitem[DPZ98]{MR1638139}
Sergio Doplicher, Claudia Pinzari, and Rita Zuccante, \emph{The
  {$C^\ast$}-algebra of a {H}ilbert bimodule}, Boll. Unione Mat. Ital. Sez. B
  Artic. Ric. Mat. (8) \textbf{1} (1998), no.~2, 263--281. \MR{1638139}

\bibitem[DR89]{MR1010160}
Sergio Doplicher and John~E. Roberts, \emph{A new duality theory for compact
  groups}, Invent. Math. \textbf{98} (1989), no.~1, 157--218. \MR{1010160}

\bibitem[Exe03a]{MR1966826}
Ruy Exel, \emph{Crossed-products by finite index endomorphisms and {KMS}
  states}, J. Funct. Anal. \textbf{199} (2003), no.~1, 153--188. \MR{1966826}

\bibitem[Exe03b]{MR2032486}
\bysame, \emph{A new look at the crossed-product of a {$C^*$}-algebra by an
  endomorphism}, Ergodic Theory Dynam. Systems \textbf{23} (2003), no.~6,
  1733--1750. \MR{2032486}

\bibitem[HP17]{MR3624399}
Michael Hartglass and David Penneys, \emph{{${\rm C}^*$}-algebras from planar
  algebras {I}: {C}anonical {${\rm C}^*$}-algebras associated to a planar
  algebra}, Trans. Amer. Math. Soc. \textbf{369} (2017), no.~6, 3977--4019.
  \MR{3624399}

\bibitem[HPN23]{2023arXiv230505072H}
Roberto Hern{\'a}ndez~Palomares and Brent {Nelson}, \emph{{Discrete Inclusions
  of $\rm{C}^{*}$-algebras}}, arXiv e-prints (2023), arXiv:2305.05072.

\bibitem[HPN24]{2024arXiv240918161H}
\bysame, \emph{{Remarks on C*-discrete inclusions}}, arXiv e-prints (2024),
  arXiv:2409.18161.

\bibitem[Jon83]{MR696688}
V.~F.~R. Jones, \emph{Index for subfactors}, Invent. Math. \textbf{72} (1983),
  no.~1, 1--25. \MR{696688}

\bibitem[Kat03]{MR2029622}
Takeshi Katsura, \emph{A construction of {$C^*$}-algebras from
  {$C^*$}-correspondences}, Advances in quantum dynamics ({S}outh {H}adley,
  {MA}, 2002), Contemp. Math., vol. 335, Amer. Math. Soc., Providence, RI,
  2003, pp.~173--182. \MR{2029622}

\bibitem[Kat04]{MR2102572}
\bysame, \emph{On {$C^*$}-algebras associated with {$C^*$}-correspondences}, J.
  Funct. Anal. \textbf{217} (2004), no.~2, 366--401. \MR{2102572}

\bibitem[KPW98]{MR1658088}
Tsuyoshi Kajiwara, Claudia Pinzari, and Yasuo Watatani, \emph{Ideal structure
  and simplicity of the {$C^\ast$}-algebras generated by {H}ilbert bimodules},
  J. Funct. Anal. \textbf{159} (1998), no.~2, 295--322. \MR{1658088}

\bibitem[KW00]{MR1624182}
Tsuyoshi Kajiwara and Yasuo Watatani, \emph{Jones index theory by {H}ilbert
  {$C^*$}-bimodules and {$K$}-theory}, Trans. Amer. Math. Soc. \textbf{352}
  (2000), no.~8, 3429--3472. \MR{1624182}

\bibitem[Pim97]{MR1426840}
Michael~V. Pimsner, \emph{A class of {$C^*$}-algebras generalizing both
  {C}untz-{K}rieger algebras and crossed products by {${\bf Z}$}}, Free
  probability theory ({W}aterloo, {ON}, 1995), Fields Inst. Commun., vol.~12,
  Amer. Math. Soc., Providence, RI, 1997, pp.~189--212. \MR{1426840}

\bibitem[Rae05]{Rae05}
Iain Raeburn, \emph{Graph algebras}, 103 of Regional conference series in
  mathematics, ISSN 0160-7642, American Mathematical Soc., 2005.

\bibitem[R{\o}r21]{2021arXiv210511899Raddendum}
Mikael R{\o}rdam, \emph{{Note on Example 5.11, Irreducible inclusions of simple
  C$^*$-algebras}}, \url{http://web.math.ku.dk/~rordam/manus/Example-5.11.pdf},
  2021.

\bibitem[R{\o}r23]{MR4599249}
\bysame, \emph{Irreducible inclusions of simple {$C^*$}-algebras}, Enseign.
  Math. \textbf{69} (2023), no.~3-4, 275--314. \MR{4599249}

\end{thebibliography}

\end{document}